\documentclass[11pt]{amsart}
\usepackage{latexsym,amssymb,amsmath,amscd,amsthm}
\numberwithin{equation}{section}
\theoremstyle{remark}
\newtheorem{theorem}{{\bf THEOREM}}[section]

\newtheorem{proposition}{{\bf PROPOSITION}}[section]
\newcommand{\bq}{\begin{equation}}
\newcommand{\bea}{\begin{array}}
\newcommand{\eea}{\end{array}}

\newcommand{\mf}{\mathfrak}

\newcommand{\ci}{\circ}

\newcommand{\gG}{\Gamma}
\newcommand{\gt}{\theta}
\newcommand{\gs}{\sigma}

\newcommand{\gd}{\delta}
\newcommand{\pp}{\partial}

\newcommand{\tl}{\tilde}

\newcommand{\bs}{\blacksquare}

\newcommand{{\DDD}}{D\!\!\!\!\!\!-}
\title{REMARKS ON Q-VIRASORO AND QKDV}
\author{Robert Carroll\\University of Illinois, Urbana, IL 61801}
\date{January, 2003\thanks{email: rcarroll@math.uiuc.edu}}
\dedicatory
{cosi lo rimembrar del dolce riso la mente mia da me medesmo scema}

\begin{document}

\bibliographystyle{plain}

\maketitle

%\newpage
%\tableofcontents
%\addcontentsline{toc}{section}{subsection}
%\setcounter{tocdepth}{3}

\section{INTRODUCTION}
\renewcommand{\theequation}{1.\arabic{equation}}
\setcounter{equation}{0}

In recent publications \cite{c2,c5,c13} we have explored derivations of
various qKdV type and noncommutative KdV (NCKdV) equations from several points of
view (e.g. differential calculi, quantum groups, hierarchies, Moyal products,
Maurer-Cartan equations, zero curvature, etc.).  Motivated in part by \cite
{k1,k2,k3} we exhibit here a derivation from q-Virasoro ideas (see e.g.
\cite{c1,c9,c12,f5,k20,k4,k5,k7,l1,s3}).  Adapting a formulation of \cite{l1}
we introduce a central term based on $\pp_q^3$ and produce 
derivations of qKdV type equations from Euler equations modeled on
the classical situation in \cite{a4,a1,c6,c20,s1,w1} (cf. also \cite{k1,m2,m1}).
One such class of equations has certain features similar to what one might imagine as
characteristic of the classical qKDV hierarchy equation (see the development in
\cite{c2} and cf. also \cite
{a2,a3,c1,c2,c4,c5,c13,f1,h1,i6,k3,t2}) and we will embed all this in a larger
framework later (see \cite{c25}).

\section{BACKGROUND}
\renewcommand{\theequation}{2.\arabic{equation}}
\setcounter{equation}{0}

First, following \cite{a1}, let $Vec(S^1)$ denote the Lie algebra of smooth
vector fields on $S^1$ and then the Virasoro algebra is $Vir=Vec(S^1)\oplus
{\bf R}={\mf W}\oplus {\bf R}$ with (note the minus sign convention involving
$f'g-fg'$)
\bq\label{2.1}
[(f(x)\pp_x,a),(g(x)\pp_x,b)]=\left((f'g-fg')\pp_x,\int_{S^1}f'g''dx\right)
\end{equation}
(${\mf W}\sim$ Witt algebra).
Here $\int_{S^1}f'g''dx$ is called the Gelfand-Fuks cocycle, where a cocycle on
a Lie algebra ${\mf g}$ is a bilinear skew symmetric form $c(\cdot,\cdot)$
satisfying ${\bf (A1)}\,\,\sum c([f,g],h)=0$ where the sum is over cyclic
permutations of $f,g,h$.  This means that $\hat{{\mf g}}={\mf g}\oplus{\bf R}$ (central
extension) with commutator $[(f,a),(g,b)]=([f,g],c(f,g))$ satisfies the Jacobi identity of
a Lie algebra.  Next one defines the Virasoro group as the set of pairs 
$(\phi(x),a)\in Diff(S^1)\oplus {\bf R}$ with multiplication law
\bq\label{2.2}
(\phi(x),a)\ci (\psi(x),b)=\left(\phi(\psi(x)),a+b+\int_{S^1}log(\phi\ci\psi(x))'
dlog\psi'(x)\right)
\end{equation}
This can be equipped with a right invariant Riemannian metric via an energy like
quadratic form on $Vir\,\,(\sim$ tangent space at the group identity) of the form
\bq\label{2.3}
H(f(x)\pp_x,a)=\frac{1}{2}\left(\int_{S^1}f^2(x)dx+a^2\right)
\end{equation}
Then the KdV equation on the circle is the evolution equation ${\bf (A2)}\,\,
\pp_tu+uu'+u'''=0$ where $'\sim\pp_x$.  More precisely the Euler equation
corresponding to geodesic flow is a 1-parameter family of KdV equations.  To see
how this arises consider ${\bf (A3)}\,\,Vir^*=\{(u(x)dx^2,c);\,\,u$ smooth on 
$S^1$ and $c\in{\bf R}\}$.  
Then
\bq\label{2.4}
<(v(x)\pp_x,a),(u(x)dx^2,c)>=\int_{S^1}v(x)u(x)dx+ac
\end{equation}
The coadjoint action of $(f\pp_x,a)\in Vir$ on $(udx^2,c)\in Vir^*$ is
($ad_v^*:\,{\mf g}^*\to {\mf g}^*,\,\,ad^*_vw(u)=w(ad_vu)$)
\bq\label{2.5}
ad^*_{(f\pp_x,a)}(udx^2,c)=(2f'u+fu'+cf''')dx^2,0)
\end{equation}
which arises from the identity
\bq\label{2.6}
<[(f\pp_x,a),(g\pp_x,b)],(udx^2,c)>=<(g\pp_x,b),ad^*_{(f\pp_x,a)}(udx^2,c)>
\end{equation}
Note here from \eqref{2.1} and \eqref{2.4}
\bq\label{2.7}
<\left(f'g-fg')\pp_x,\int_{S^1}f'g''dx\right),(udx^2,c)>=\int_{S^1}(f'g-fg')udx+c\int_{S^1}f'g''dx
\end{equation}
while from \eqref{2.4} and \eqref{2.5}
\bq\label{2.8}
<(g\pp_x,b),ad^*_{f\pp_x,a)}(udx^2,c)>=<(g\pp_x,b),((2f'u+fu'+cf''')dx^2,0)>
\end{equation}
Now for $S^1$ there are no boundary terms in integration (for single valued
functions) so, integrating by parts,
\bq\label{2.9}
\int_{S^1}g(2f'u+fu'+cf''')dx=\int_{S^1}[u(gf'-fg')+cf'g'']dx
\end{equation}
in agreement with \eqref{2.7}.
Now a function H on ${\mf g}=Vir$ determines a tautological inertia operator 
$A:\,Vir\to Vir^*:\,(u\pp_x,c)\to (udx^2,c)$ and hence a quadratic Hamiltonian
on $Vir^*$ via
\bq\label{2.10}
H(udx^2,c)=\frac{1}{2}\left(\int_{S^1}udx^2+c^2\right)=\frac{1}{2}<(u\pp_x,c),(udx^2,c)>=
\frac{1}{2}<(u\pp_x,c),A(u\pp_x,c)>
\end{equation}
Following \cite{k1} the corresponding Euler equation is
$\dot{m}=-ad^*_{A^{-1}m}m\,\,(m\in\hat{{\mf g}}$) which here takes the form
\bq\label{2.11}
\pp_t(udx^2,c)=-ad^*_{A^{-1}(udx^2,c)}(udx^2,c)
\end{equation}
which becomes via \eqref{2.5} with ${\bf
(A4)}\,\,(f\pp_x,a)=A^{-1}(udx^2,c)=(u\pp_x,c)$
\bq\label{2.12}
\pp_tu=-2u'u-uu'-cu'''=-3uu'-cu'''
\end{equation}
where $c$ is independent of time.  One notes also that without the central
extension of $Vec(S^1)$ we arrive in the same manner at a nonviscous Burger's
equation ${\bf (A5)}\,\,\pp_tu=-3uu'$.

\section{Q-VIRASORO}
\renewcommand{\theequation}{3.\arabic{equation}}
\setcounter{equation}{0}

For q-Virasoro we go directly to \cite{l1} and refer to the bibliography in
the references given above in the introduction for other approaches.  Thus
work on $S^1$ with ($q\ne 0,\pm 1$)
\bq\label{3.1}
\pp_qz=\frac{q^mz^m-q^{-m}z^m}{(q-q^{-1}z}=z^{m-1}[m];\,\,[m]=
\frac{q^m-q^{-m}}{q-q^{-1}}
\end{equation}
We adapt the formalism of \cite{l1} as follows.  Let $D_n=-z^{n+1}\pp$ with
$\pp:\,z^m\to q^m[m]z^{m-1}$ so $\pp\sim\pp_q\tau$ where $\tau f(z)=f(qz)$.
Generally we will think of $z=e^{i\gt}\in S^1$ so $(1/2\pi
i)\int_{S^1}z^ndz=(1/2\pi)\int z^{n+1}d\gt=\gd_{(-1,0)}$ which will be written
as ${\bf (A6)}\,\,\int z^n=\gd_{(-1,0)}$.  Write also ${\bf (A7)}\,\,
\ell_n\sim z^{n+1}\pp=-D_n=z^{n+1}\pp_q\tau$.  It is known that q-brackets
are needed now where
\bq\label{3.2}
[\ell_m,\ell_n]_q=q^{m-n}\ell_m\ell_n-q^{n-m}\ell_n\ell_m=[m-n]\ell_{m+n}
\end{equation}
For a central term in a putative $Vir_q$ one wants (cf. \cite{a5,c14,l1})
a formula ${\bf (A8)}\,\,c[m+1][m][m-1]\gd_{m+m,0}$ (see Remark 3.1 for an optimal term).
First we want to
formulate the q-bracket in terms of vector fields as follows (the central
term will be added later in a somewhat ad hoc manner).  This can be done as a
direct calculation (note $\pp_qf=(\tau f)\pp_q+(\pp_qf)\tau^{-1}$)
\bq\label{3.3}
[z^n\pp,z^m\pp]_q\sim q^{n-m}z^n\pp(z^m\pp)-q^{m-n}z^m\pp(z^n\pp)=
q^{n-m}z^nq^m\pp_qz^m\tau\pp-
\end{equation}
$$-q^{m-n}z^mq^n\pp_qz^n\tau\pp=q^{n-m}z^nq^m(q^mz^m\pp_q+[m]_qz^{m-1}\tau^{-1})\tau\pp
-q^{m-n}z^mq^n(q^nz^n\pp_q+$$
$$+[n]_qz^{n-1}\tau^{-1})\tau\pp=(q^n[m]-q^m[n])z^{m+n-1}\pp=
[n-m]z^{m+n-1}\pp$$
Let now $v\sim \sum a_nz^n$ and $w\sim \sum b_mz^m$; then we define a bracket in
$Vec(S^1)$ via
\bq\label{3.4}
[v\pp,w\pp]_q\sim
(\tau v)\pp(\tau^{-1}w)\pp-(\tau w)\pp(\tau^{-1}v)\pp=(\tau v)\pp_q\tau(\tau^{-1}w)\pp-
\end{equation}
$$-(\tau w)\pp_q\tau(\tau^{-1}v)\pp=(\tau v)(\pp_qw\tau\pp)-(\tau w)(\pp_qv\tau\pp)=$$
$$=\sum a_nq^nz^n(\sum b_m(q^mz^m\pp_q+[m]z^{m-1}\tau^{-1}))\tau\pp-\sum b_mz^mq^m
(\sum a_n(q^nz^n\pp_q+$$
$$+[n]z^{n-1}\tau^{-1}))\tau\pp=\sum a_nb_m(q^n[m]-q^m[n])z^{m+n-1})\pp=\sum a_nb_m
[n-m]z^{m+n-1}\pp$$
\begin{proposition}
From \eqref{3.4} we have a correspondence
\bq\label{3.5}
v'w-vw'\sim -[v\pp_x,w\pp_x]\sim -[v\pp,w\pp]_q=-\{(\tau v)(\pp_qw)-(\tau
w)(\pp_qv)\}\tau
\end{equation}
\end{proposition}
\indent
{\bf REMARK 3.1.}
In \cite{l1} one defines the q-analogue of the enveloping algebra 
of the Witt algebra ${\mf W}$ as the associative algebra ${\mf U}_q({\mf W})$
having generators $\ell_m\,\,(m\in{\bf Z})$ and relations \eqref{3.2}.  The q-deformed
Virasoro algebra is defined as the associative algebra ${\mf U}_q(Vir)$ having generators
$\ell_m\,\,(m\in {\bf Z})$ and relations ($q\ne$ root of unity)
\bq\label{alpha}
q^{m-n}\ell_m\ell_n-q^{n-m}\ell_n\ell_m=[m-n]\ell_{m+n}+\gd_{m+n,0}\frac
{[m+1][m][m-1]}{[2][3]<m>}\hat{c}
\end{equation}
where $<m>=q^m+q^{-m}$ and $\hat{c}\ell_m=q^{2m}\ell_m\hat{c}$ (thus $\hat{c}$ is an
operator which we examine below and we refer to \cite{a5,c14,l1} for the central term).  
Then
${\mf U}_q(Vir)\sim Vir_q$ is a {\bf Z} graded algebra with $deg(\ell_m)=m$ and
$deg(\hat{c})=0$.  One also introduces in 
\cite{l1} a larger algebra ${\mf U}(V_q)=$ associative algebra generated by $J^{\pm 1},
\,\hat{c},\,d_m\,\,(m\in{\bf Z})$ with relations
\bq\label{beta}
JJ^{-1}=J^{-1}J=1;\,\,Jd_mJ^{-1}=q^md_m;\,\,\hat{c}J=J\hat{c};\,\,\hat{c}d_m=q^md_m\hat{c};
\end{equation}
$$q^md_md_nJ-q^nd_nd_mJ=[m-n]d_{m+n}+\gd_{m+n,0}\frac{[m+1][m][m-1]}{[2][3]<m>}\hat{c}$$
The subalgebra of ${\mf U}(V_q)$ generated by $\ell_m'=d_mJ$ and
$\hat{c}'=\hat{c}J\,\,(m\in {\bf Z})$ is the same as ${\mf U}_q(Vir)$.
We will treat \eqref{alpha} as displaying the natural (preferred) central term for our
purposes.
$\hfill\bs$

\section{CALCULATIONS}
\renewcommand{\theequation}{4.\arabic{equation}}
\setcounter{equation}{0}

Now we mimic the framework of Section 2 and it is interesting to note that an
ordinary integral $\int_{S^1}$ will suffice.  One does not need a Jackson type
integral in order to deal with integration by parts.  Thus we observe that
\bq\label{4.1}
\int_{S^1}f=\int\sum f_nz^n=f_{-1};\,\,\int\pp_qf=\frac{1}{q-q^{-1}}
\int\frac{f(qz)-f(q^{-1}z)}{z}=
\end{equation}
$$=\frac{1}{q-q^{-1}}\int \sum f_nz^{n-1}(q^n-q^{-n})=\frac{1}{q-q^{-1}}
(f_0-f_0)=0$$
Since $\pp_q(fh)=(\tau f)(\pp_qh)+(\pp_qf)(\tau^{-1}h)$ we have an
integration by parts formula
\bq\label{4.2}
\int(\tau f)(\pp_qh)=-\int (\pp_qf)(\tau^{-1}h)\Rightarrow \int f\pp_q(\tau
h)=-\int \pp_q(\tau^{-1}f)h
\end{equation}
This can be written as (recall $\pp\sim \pp_q\tau$) ${\bf (A9)}\,\,\int f\pp h
=-\int h\hat{\pp}f$ for $\hat{\pp}=\pp_q\tau^{-1}$.  Now we think of
${\mf U}_q(Vir)$ with elements $(f\pp,a)$ as in \eqref{2.1} with
\bq\label{4.3}
[(f\pp,a),(g\pp,b)]=(-[f\pp,g\pp]_q\pp,\int (\tau\pp^3f)(\tau g)\hat{c})
\end{equation}
The central term is defined tentatively via ${\bf (A10)}\,\,\int
(\tau\pp^3f)(\tau g)\hat{c}=
\psi(f\pp,g\pp)$ where one has ${\bf (A11)}\,\,\psi(f\pp,g\pp)=
q^{-1}\int g\pp^3f\hat{c}$ since
\bq\label{4.4}
\int \tau a\tau b=\int \sum a_nb_mq^{n+m}z^{n+m}=\sum a_nb_{-n-1}q^{-1}=
q^{-1}\int ab
\end{equation}
We will want to put the central operator $\hat{c}$ into the integral 
{\bf (A10)} or {\bf (A11)}, acting on f, and will see below that $\hat{c}\sim \tau^2$ for
example and $\tau^2F(z)= F(q^2z)\tau^2$ so it eventually automatically passes to the
right in our qKdV type equations.  Hence for the moment think of $\hat{c}=\tau^2$ put
into \eqref{4.3} or {\bf (A10)} via e.g. $\tau\pp^3\tau^{-2}\hat{c}f\equiv \tau\pp^3f$
and ignored at the end except when exhibiting formulas like \eqref{alpha} on generators
(see also Remark 4.5).
\\[3mm]\indent
Now duality as in \eqref{2.4} will be expressed here via ${\bf (A12)}\,\,
<(v\pp,a),(u,c)>=\int v\tau^{-1}u+ac$ and hence from \eqref{3.5} 
\bq\label{4.5}
<[(f\pp,a),(g\pp,b)],(u,c)>=\int
-[f\pp,g\pp]_q\tau^{-1}u+c\psi(f\pp,g\pp)=
\end{equation}
$$=\int (\tau g)(\pp_qf)u-(\tau
f)(\pp_qg)u+c\int (\tau\pp^3f)(\tau g)\hat{c}$$
This puts us in the framework of \eqref{2.7}, \eqref{2.9}, etc.  Finally note ${\bf
(A13)}\,\,\pp_q(g(\tau^2f)\tau u)=(\pp_qg)(\tau f)u+(\tau g)\pp_q((\tau^2f)\tau u)$ so
\eqref{4.5} becomes 
\bq\label{4.6}
\int (\tau g)[(\pp_qf)u+\pp_q((\tau^2f)\tau u)+c\tau\pp^3f]
\end{equation}
\begin{theorem}
In the spirit of Section 2 \eqref{4.6} leads to a tentative qKdV type equation for 
$f=u$ (note $\pp_q\tau =q\tau\pp_q$)
\bq\label{4.7}
u_t=-c\tau\pp^3u-(\pp_qu)u-\pp_q((\tau^2u)\tau u)=-c\tau\pp^3u-(1+q\tau)^2
u\pp_qu
\end{equation}
where $\pp^3u=\pp_q\tau\pp_q\tau\pp_q\tau u$ (cf. also remarks and results below
especially where the operator $\hat{c}$ is put into the equation more meaningfully
and the central term follows \eqref{alpha} acting on generators).
\end{theorem}
\indent
{\bf REMARK 4.1.}
$\psi(f\pp,g\pp)$ appears to be a perfectly satisfactory central term even though it
does not seem to be a cocycle on two counts
and thus the theorem uses a reduced structure
for its derivation (e.g. there will not be a Jacobi identity with
brackets \eqref{4.3} - for more on this see Remarks 4.2 - 4.6). 
First via
${\bf (A11)}\,\,\psi(f\pp,g\pp)=q^{-1}\int g\pp^3f=q^{-1}\hat{\psi}
(f\pp,g\pp)$ and an 
elementary calculation gives ${\bf (A14)}\,\,\int g\pp^3f=-\int
f\hat{\pp}^3g$ with $\hat{\pp}$ as in {\bf (A9)}.  Further for $f=\sum f_{n+1}
z^{n+1}$ and $g=\sum g_{m+1}z^{m+1}$
\bq\label{4.8}
\int g\pp^3f=\int \sum f_{n+1}g_{m+1}[n+1][n][n-1]q^{3n}x^{n+m-1}=
\sum f_{n+1}g_{-n+1}[n+1][n][n-1]q^{3n};
\end{equation}
$$\int f\hat{\pp}^3g=\int\sum f_{n+1}g_{m+1}[m+1][m][m-1]q^{-3m}x^{m+n-1}=$$
$$=\sum f_{n+1}g_{-n+1}[-n+1][-n][-n-1]q^{3n}=-\int g\pp^3f$$
However the antisymmetry condition $\psi(f\pp,g\pp)=-\psi(g\pp,f\pp)$ does not hold since
\bq\label{star}
\int f\pp^3g=\sum f_{n+1}g_{-n+1}[-n+1][-n][-n-1]q^{-3n}
\end{equation}
We note however that $\int (\pp^3z^{n+1})z^{m+1}=
[n+1][n][n-1]\int z^{m+n-1}=[n+1][n][n-1]\gd_{m+n,0}$ as in {\bf (A8)}
but this falls short of \eqref{alpha} by an n dependent term $<n>^{-1}$.  The remaining
cocycle condition can be written as
\bq\label{4.9}
\hat{\psi}(f\pp,[g\pp,h\pp]_q\pp)+\hat{\psi}(h\pp,[f\pp,g\pp]_q\pp)
+\hat{\psi}(g\pp,[h\pp,f\pp]_q\pp)=0
\end{equation}
and one can write out the terms 
directly to see that cancellation does not occur.$\hfill\bs$
\\[3mm]\indent
{\bf REMARK 4.2.}
One can salvage a bit here by taking as central term
\bq\label{4.16}
\tl{\psi}(f\pp,g\pp)=\int (\tau\pp^3\tau^{-3}f)(\tau
g)=q^{-1}\int(\pp^3\tau^{-3}f)g=
\end{equation}
$$=q^{-4}\int\sum f_{n+1}g_{m+1}[n+1][n][n-1]x^{n+m-1}=q^{-4}\sum f_{n+1}g_{-n+1}
[n+1][n][n-1]$$
For this we have
\bq\label{4.17}
\tl{\psi}(g\pp,f\pp)=q^{-4}\sum f_{n+1}g_{-n+1}[-n+1][-n][-n-1]=-\tl{\psi}(f\pp,
g\pp)
\end{equation}
so antisymmetry is realized.  However a little calculation shows that the cocycle
condition \eqref{4.9} again does not hold.
$\hfill\bs$
\\[3mm]\indent
Thus let us indicate what algebraic structure is possesed by $Vec_q(S^1)\oplus{\bf R}$ in
our constructions with the bracket of Proposition 3.1 and the central term added as
in \eqref{4.16} since this is at least antisymmetric.  First one states
\begin{theorem}
Using $\tl{\psi}$ as in Remark 4.2 one gets a tentative qKdV equation 
\bq\label{4.24}
u_t=-c\tau\pp^3\tau^{-3}u-(1+q\tau)^2u\pp_qu
\end{equation}
\end{theorem}
\indent
{\bf REMARK 4.3.}
Regarding structure one does not have a central extension (as with $Vir=Vec(S^1)\oplus
{\bf R}$ in the classical case) but there is an antisymmetric bracket (cf. \eqref{4.3})
\bq\label{4.25}
[(g\pp,b),(f\pp,a)]=(-[g\pp,f\pp]_q,\tl{\psi}(g\pp,f\pp))=-[(f\pp,a),(g\pp,b)]=
\end{equation}
$$=-(-[f\pp,g\pp]_q,\tl{\psi}(f\pp,g\pp))=-([g\pp,f\pp]_q,\tl{\psi}(f\pp,g\pp))$$
However without a cocycle we do not have a Jacobi identity.  Nevertheless a dual
structure can be defined as in {\bf (A12)} and manipulated as in Section 2 (with no need
to refer to quadratic differentials, etc.).  Further the structure imposed by $\tl{\psi}$
does give (reinserting $\hat{c}$ on generators as discussed above after \eqref{4.4})
\bq\label{4.26}
\tl{\psi}(\ell_n,\ell_m)=\tl{\psi}(z^{n+1}\pp,z^{m+1}\pp)=q^{-4}[n+1][n][n-1]
\gd_{m+n,0}\hat{c}
\end{equation}
as stipulated in {\bf (A8)} so we are speaking of 
$Vec_q(S^1)\oplus{\bf R}$ defined on generators without specifying (or needing) any
additional structure.$\hfill\bs$
\\[3mm]\indent
{\bf REMARK 4.4.}
We add a few facts designed to clarify the amount of structure used in Section 2
and needed in the q-theory with central term arising via $\tl{\psi}$ for example.
Thus in Section 2 $Vir=Vec(S^1)\oplus {\bf R}={\mf W}\oplus {\bf R}\sim\hat{{\mf W}}$
via \eqref{2.1} where the cocycle can be written as $\int_{S^1}f'''g$ equally well.
In this context $Vir\sim\hat{{\mf W}}$ (standard central extension).  One recalls that
the Gelfand-Fuks cocycle (and hence the Virasoro algebra) measures the deformation
of the projective structure on $S^1={\bf R}P^1$ by diffeomorphisms (cf. \cite{a1})
and is thus a rather important and specific kind of cocycle.  One could therefore
argue that KdV is a very important equation with significant geometric content
(the importance would appear to be well known but the geometric significance for
qKdV remains to be expressed in terms of quantum groups). 
Now consider the
q-theory defined with central term arising via $\tl{\psi}$.  For lack of a
cocycle we don't have a central extension.  However
the definition {\bf (A12)} of a $(Vec(S^1)\oplus{\bf R})^*$ (vector space
dual) related to {\bf (A3)} is straightforward as is the generalization of $ad^*$ in
\eqref{2.6}. No algebraic structure is used or needed in tracing this development as in
Section 4 although the significance of the resulting equations seems dependent on the
form of central term and would presumably be enhanced in the presence of a ``q-cocycle"
of some sort.
$\hfill\bs$
\\[3mm]\indent
{\bf REMARK 4.5.}
Now from \cite{l1} we know that ${\mf U}_q({\mf W})$ is
an associative algebra with generators $\ell_m$ and q-bracket as in \eqref{3.2}.  The
central term in 
\cite{l1} is ${\bf (A15)}\,\,([m+1][m][m-1]/[2][3]<m>)\gd_{m+n,0}c$ where
$<m>=q^m+q^{-m}$.
In our situation we have in \eqref{4.26} a term $q^{-4}[n+1][n][n-1]\hat{c}\sim -q^{-4}
[m+1][m][m-1]\hat{c}$ and to bring this into line one takes first $\hat{c}=
\tau^2$ where ${\bf (A16)}\,\,\tau^2\ell_m=\tau^2(z^{m+1}\pp_q\tau)=q^{2m}z^{m+1}\pp_q
\tau\tau^2=q^{2m}\ell_m\tau^2$ (note $\pp_q\tau=q\tau\pp_q$ and
$\pp_q\tau^{-1}=q^{-1}\tau^{-1}\pp_q$).  Then to get a factor 
$1/<m>$ note simply ${\bf (A17)}\,\,
(\tau +\tau^{-1})z^m=<m>z^m$ and 
consequently
\bq\label{delta}
(\tau+\tau^{-1})^{-1}z^m=<m>^{-1}z^m
\end{equation}
So let us use $\hat{c}\sim \tau^2$ to get $\hat{c}\ell_m=q^{2m}\ell_m\hat{c}$ as
indicated in Remark 3.1.  Then we will work $(\tau+\tau^{-1})^{-1}$ into the calculation
to obtain $<m>^{-1}$.  Finally pick a $c=c'q^4/[2][3]$ for arbitrary $c'$.  Following
up as remarked after \eqref{4.4} we note
that the presence of an operator $\hat{c}\sim \tau^2$ would complicate our putative qKdV
equations such as \eqref{4.7} or \eqref{4.24} if it is left hanging on the end.
On the other hand we don't want to leave it out of the calculation so we move it to the
left and let it work on $z^{n+1}$ for example or on f.  Thus we set 
${\bf (A18)}\,\,\hat{c}(z^{n+1})=q^{2n+2}z^{n+1}$ (omitting the trailing $\hat{c}$
except when needed) and this is illustrated below.
$\hfill\bs$
\\[3mm]\indent
{\bf REMARK	 4.6.}
Consider a generic situation and define now a new $\psi$ via ($\hat{c}\sim \tau^2,\,\,
\pp=\pp_q\tau$)
\bq\label{4.29}
\psi(\ell_n,\ell_m)=q^6\int (\tau z^{m+1})\tau(\pp^2\hat{c}(\tau+\tau^{-1})^{-1}\pp
\tau^{-5}z^{n+1})=
\end{equation}
$$=q^5\int z^{m+1}(\pp^2\hat{c}(\tau+\tau^{-1})^{-1}\pp_q\tau^{-4}z^{n+1})=q
\int q^{-2n}z^{m+1}(\pp^2z^n)\frac{[n+1]}{<n>}=$$
$$=q\int q^{-2n}z^{m+1}q^n\pp z^{n-1}\frac{[n+1][n]}{<n>}=\frac{[n+1][n][n-1]}{<n>}
\gd_{m+n,0}\hat{c}$$
where $\hat{c}=\tau^2$ is restored at the right in the last equation (since $\tau
F(x)=F(qx)\tau$) so one recovers \eqref{alpha} up to a factor of $[2][3]$ (which
can be recovered by writing $\hat{c}=c'\tau^2$ for a suitable constant $c'$,
independent of n).  Applying
this more generally we write for $\hat{c}=(q^{-6}/[2][3])c'\tau^2$ ($c'$ arbitrary)
\bq\label{4.30}
\psi(f\pp,g\pp)=\int (\tau g)\tau(\pp^2\hat{c}(\tau+\tau^{-1})^{-1}\pp\tau^{-5}f)
\end{equation}
This leads to
\begin{theorem}
A possibly canonical qKdV type equation can be obtained in the form
\bq\label{4.31}
u_t=-c''\tau\pp^2\tau^2(\tau+\tau^{-1})^{-1}\pp\tau^{-5}u-(1+q\tau)^2u\pp_qu
\end{equation}
where $\pp=\pp_q\tau$.
\end{theorem}
\indent
{\bf REMARK 4.7.}
It would be possible to vary the third derivative term in various ways by repositioning
of $\tau^2$ and $(\tau+\tau^{-1})^{-1}$ and we have chosen this way for convenience. 
One very interesting feature of this equation involves the $(\tau+\tau^{-1})^{-1}$ term.
This serves to bring the equation into line with the structure of the classical
qKdV hierarchy equation suggested implicitly in \cite{c2}, since expansion of
$(\tau+\tau^{-1})^{-1}$ would involve an infinite number of terms.  This will be
investigated further (cf. also \cite{a2,a3,c5,c13,f1,h1,i6,k3}); a few
further details are provided in Remarks 4.9, 4.11, and 4.12.$\hfill\bs$
\\[3mm]\indent
{\bf REMARK 4.8.}
We check now some properties of $\psi$ as defined in \eqref{4.30}.  Thus write
\bq\label{4.32}
\psi(f\pp,g\pp)=q^{-1}c''\int g\pp^2\tau^2(\tau+\tau^{-1})^{-1}\pp\tau^{-5}f=
\end{equation}
$$=q^{-1}c''\int\sum g_{m+1}z^{m+1}f_{n+1}\pp^2\tau^2(\tau+\tau^{-1})^{-1}\pp_q\tau^{-4}
z^{n+1}=$$
$$=q^{-1}c''\int\sum g_{m+1}f_{n+1}z^{m+1}\pp^2\tau^2(\tau+\tau^{-1})^{-1}q^{-4n-4}
[n+1]z^n=$$
$$=q^{-1}c''\int\sum g_{m+1}f_{n+1}z^{m+1}q^{-4n-4}\frac{[n+1]}{<n>}\pp^2q^{2n}z^n=$$
$$=q^{-1}c''\int\sum
g_{m+1}f_{n+1}z^{m+1}q^{-2n-4}\frac{[n+1][n]}{<n>}\pp_q\tau q^nz^{n-1}=$$
$$=q^{-5}c''\int\sum g_{m+1}f_{n+1}
q^{-2n}z^{m+1}q^n\frac{[n+1][n][n-1]}{<n>}q^{n-1}z^{n-2}=$$
$$=q^{-6}\int\sum g_{m+1}f_{n+1}z^{m+n-1}\frac{[n+1][n][n-1]}{<n>}=q^{-6}c''
\sum f_{n+1}g_{-n+1}\Xi(n)$$
On the other hand ($n\to -n$)
\bq\label{4.33}
\psi(g\pp,f\pp)=q^{-6}c''\sum g_{n+1}f_{-n+1}\Xi(n)=-q^{-6}c''\sum f_{n+1}g_{-n+1}
\Xi(n)
\end{equation}
since $\Xi(-n)=-\Xi(n)$.  Hence $\psi$ is antisymmetric.  As for the cocycle condition
\eqref{4.9} we consider e.g. (cf. \eqref{3.4} and \eqref{4.32})
\bq\label{4.34}
\psi(f\pp,[g\pp,h\pp]_q)=q^{-6}c''\int\sum f_{n+1}g_{m+1}h_{s+1}[s-m]\frac{[n+1][n][n-1]}
{<n>}z^{m+n+s-1}
\end{equation}
Similarly one will have terms
\bq\label{4.35}
\psi(h\pp,[f\pp,g\pp]_q)=q^{-6}c''\int\sum f_{n+1}g_{m+1}h_{s+1}\frac{[s+1][s][s-1]}{<s>}
z^{m+n+s-1}[m-n]
\end{equation}
\bq\label{4.36}
\psi(g\pp,[h\pp,f\pp]_q)=q^{-6}c''\int\sum f_{n+1}g_{m+1}h_{s+1}[n-s]\frac{[m+1][m][m-1]}
{<m>}z^{m+n+s-1}
\end{equation}
Thus one would need here
\bq\label{4.37}
\Xi(n)[s-m]+\Xi(s)[m-n]+\Xi(m)[n-s]=0\,\,\,(n+m+s=0)
\end{equation}
We checked previously (cf. Remark 4.2) that this does not hold when the $<n>,\,\,<m>,$
and $<s>$ terms are absent and leave this open for the moment.$\hfill\bs$
\\[3mm]\indent
{\bf REMARK 4.9.}
In \cite{c2} we exhibited the qKdV hierarchy equation in the form
\bq\label{4.38}
\pp_tu=[L^3_{+},L^2]_0=[D_q^3+w_2D_q^2+w_1D_q+w_0,D_q^2+u_1D_q+u]_0=
\end{equation}
$$(D_q^3u)+w_2(D_q^2u)+w_1(D_qu)-[(D_q^2w_0)+u_1(D_qw_0)]$$
where $u_1=(q-1)xu$ and $w_i=w_i(u)$.  Here $D_qf=[f(qx)-f(x)]/(q-1)x$ and $L=D_q+
s_0+s_1D_q^{-1}+\cdots$.  However
\bq\label{4.39}
w_2=\tau^2s_0+u_1=\tau^2s_0+\tau s_0+s_0;
\end{equation}
$$w_1=(q+1)(\tau D_qs_0)+\tau^2s_1+[(\tau s_0)+s_0](\tau s_0)+
u;$$
$$w_0=D_q^2s_0
+(q+1)(\tau D_qs_1)+u_1D_qs_0+u_1(\tau s_1)+us_0+\tau^2s_2$$
and the determination of the $s_i$ requires infinite series calculations based e.g. on
formulas like
\bq\label{4.40}
s_1+\tau s_1=u-D_qs_0-s^2_0=f\Rightarrow s_1(x)=\sum_0^{\infty}(-1)^nf(q^nx)
\end{equation}
Similarly ${\bf (A19)}\,\,\tau s_2+s_2=-D_qs_1-s_0s_1-s_1\tau^{-1}s_0$, etc.  This seems
to suggest that an explicit form for a qKdV equation of this type with all coefficients
specified could involve an infinite series of the type appearing in \eqref{4.31} upon 
expansion of $(\tau+\tau^{-1})^{-1}$ (see here Remark 4.11).$\hfill\bs$
\\[3mm]\indent
{\bf REMARK 4.10.}
The use of q-brackets and formulas like \eqref{alpha} seem to bypass the quest for an
embedding of $Vec_q(S^1)$ into qPSDO involving the logarithmic cocycle of \cite{k3}.  It
would be interesting to write out explicitly a form for qKdV in the hierarchy theory
using the embedding of \cite{k3}.$\hfill\bs$
\\[3mm]\indent
{\bf REMARK 4.11.}
Following Remark 4.7 we look at an alternative rendering of \eqref{4.30} - \eqref{4.31}.
Thus consider for $\hat{c}=\tau^2$
\bq\label{4.41}
\psi(f\pp,g\pp)=\int (\tau g)(\tau\tau^{-5}\pp^2(\tau+\tau^{-1})^{-1}\pp(\tau^2c'f)
\end{equation}
from which we obtain ${\bf (A20)}\,\,\psi(\ell_n,\ell_m)=c'q^6([n+1][n][n-1]/<n>)
\gd_{m+n,0}\hat{c}$ so one chooses $c'=c''q^{-6}/[2][3]$ etc..  The corresponding 
qKdV equation would be formally
\bq\label{4.42}
u_t=-c''\tau^{-4}\pp^2(\tau+\tau^{-1})^{-1}\pp\tau^2 u-(1+q\tau)^2u\pp_qu=
\end{equation}
$$=-(1+q\tau)^2
u\pp_qu-c''\tau^{-4}\pp_q\tau\pp_q\tau\sum_0^{\infty}(-1)^n\tau^{2n+3}\pp_q\tau u$$
A similar formula would apply for \eqref{4.31}.$\hfill\bs$
\\[3mm]\indent
{\bf REMARK 4.12.}
We add a few embellishments as follows (cf. the CR note in \cite{l1}).  Define 
$\gG =(q^{-1}\tau+q\tau^{-1})=q\tau^{-1}(1+q^{-2}\tau^2)$ with $\gG^{-1}=(\tau/q)
(1+q^{-2}\tau^2)^{-1}$ so $\gG z^{n+1}=<n>z^{n+1}$.  Then consider ($\hat{c}\sim \tau^2$)
\bq\label{4.43}
\psi(f\pp,g\pp)=\int (\tau g)\tau(\tau^{-5}\pp^3c'\hat{c}\gG^{-1}f)
\end{equation}
It follows that $\psi(\ell_n,\ell_m)=(c'q^6[n+1][n][n-1]/<n>)\gd_{m+n,0}\hat{c}$ and
\bq\label{4.44}
u_t=-(1+q\tau)^2u\pp_qu-c''\tau^{-4}\pp_q\tau\pp_q\tau\pp_q\sum_0^{\infty}
(-1)^nq^{-2n-1}\tau^{2n+4}u
\end{equation}
Now it is stated in \cite{l1} (CR note) that $Vir_q$ is
the universal quantum central extension of ${\mf W}_q$ but this is not clear (note that
the $d_m$ generators of Remark 3.1 are used in this context).  For $Witt_q\sim{\mf W}_q$
however there is a quantum Jacobi identity based on 
\bq\label{4.45}
[m-n][m+n-p]<p>+[n-p][n+p-m]<m>+[p-m][p+m-n]<n>=0
\end{equation}
$\gG(\ell_p)=<p>\ell_p$.  Then writing $\gG(\ell_p)=<p>\ell_p$ and
$\gs_m=([m+1][m][m-1]/<m>[2][3])$ we have
\bq\label{4.46}
[[\ell_m,\ell_n]_q,\gG(\ell_p)]_q=[m-n][\ell_{m+n},<p>\ell_p]_q+\gs_m\gd_{m+n,0}
[\hat{c},\gG(\ell_p)]_q
\end{equation}
However ($deg(\hat{c})=0$) $[\hat{c},\ell_p]=\hat{c}\ell_pq^{-p}-\ell_p\hat{c}q^p=
\ell_p(q^p-q^p)\hat{c}=0$ so, since $\hat{c}\ell_p=q^{2p}\ell_p\hat{c}$,
\bq\label{4.47}
[[\ell_m,\ell_n]_q,\gG(\ell_p)]_q=[m-n][m+n-p]<p>\ell_{m+n+p}+\gs_{m+n}
[m-n]<p>\gd_{m+n+p,0}\hat{c}
\end{equation}
Hence (via \eqref{4.45})
\bq\label{4.48}
[[\ell_m,\ell_n]_q,\gG(\ell_p)]_q+[[\ell_n,\ell_p]_q,\gG(\ell_m)]_q+[[\ell_p,\ell_m]_q,
\gG(\ell_n)]_q=
\end{equation}
$$=\{\gs_{m+n}[m-n]<p>+\gs_{n+p}[n-p]<m>+\gs_{p+m}
[p-m]<n>\}\gs_{m+n+p}\hat{c}=0$$
The right side comes only from the Virasoro central term so we get zero for the 
${\mf W}_q$ algebra.  To check Jacobi for $Vir_q$ one looks at
\bq\label{4.49}
\gs_{m+n}[m-n]<p>+\gs_{n+p}[n-p]<m>+\gs_{p+m}[p-m]<n>
\end{equation}
for $p+m+n=0$.  This is simply
\bq\label{4.50}
[m+n+1][m+n][m+n-1]+[n+p+1][n+p][n+p-1]+[p+m+1]\times
\end{equation}
$$\times [p+m][p+m-1]=
[m+m+1][m+n][m+n-1]-
[m+1][m][m-1]-[n+1][n][n-1]$$
but this does not seem to vanish.  Possible
cocycle formulas here for $\psi$ will run into the same calculations as in
Remark 4.8.$\hfill\bs$

\end{document}